\documentclass[letterpaper, 10 pt, conference]{ieeeconf}  % Comment this line out if you need a4paper

\IEEEoverridecommandlockouts                              % This command is only needed if 
\usepackage{graphicx}
\usepackage{stfloats}
\usepackage{tabularx}

\usepackage{epsfig} % for postscript graphics files
\usepackage{epstopdf}
\usepackage{amssymb}  % assumes amsmath package installed

\usepackage{cite}
\usepackage{url}
\usepackage[hidelinks]{hyperref}
\usepackage[hyphenbreaks]{breakurl}
\usepackage{booktabs}
\usepackage[american,fulldiode]{circuitikz} %%Circuit drawing tools

\usepackage{algorithm}
\usepackage[noend]{algpseudocode}
\usepackage{float}
\newfloat{algorithm}{t}{lop}

\usepackage[caption=false,font=footnotesize]{subfig}

\usepackage{tablefootnote}

\renewcommand{\leftmark}{{Preprint No: ADP-15-35/T937}}
\renewcommand{\rightmark}{{Preprint No: ADP-15-35/T937}}

\usepackage{amsmath} % assumes amsmath package installed
\usepackage{amssymb}  % assumes amsmath package installed

\renewcommand{\Re}{\operatorname{Re}}
\renewcommand{\Im}{\operatorname{Im}}
\newcommand{\msosr}[1]{\mbox{\ensuremath{\text{MSOS}_{#1}}-$\mathbb{R}$}}
\newcommand{\msosc}[1]{\mbox{$\text{MSOS}_{#1}$-$\mathbb{C}$}}

\title{\LARGE \bf
Recent Advances in Computational Methods for the \\ Power Flow Equations
}

\author{Dhagash Mehta$^{1}$, Daniel K. Molzahn$^{2}$, and Konstantin Turitsyn$^{3}$% <-this % stops a space
\thanks{$^{1}$University of Notre Dame, Dept. of Applied and Computational Mathematics and Statistics; 
University of Adelaide, School of Physical Sciences, Dept. of Physics, Centre for the Subatomic Structure of Matter. 
Support from NSF-ECCS award ID 1509036 and an Australian Research Council DECRA fellowship no. DE140100867.
        {\tt\small dmehta@nd.edu}, \mbox{\scriptsize Preprint No: ADP-15-35/T937}}%
\thanks{$^{2}$Argonne National Laboratory, Energy Systems Division.
        {\tt\small dmolzahn@anl.gov}}%
\thanks{$^{3}$ Department of Mechanical Engineering, Massachusetts
Institute of Technology.
        {\tt\small turitsyn@mit.edu}}%
}

\begin{document}

\maketitle
\thispagestyle{empty}
\pagestyle{empty}

%%%%%%%%%%%%%%%%%%%%%%%%%%%%%%%%%%%%%%%%%%%%%%%%%%%%%%%%%%%%%%%%%%%%%%%%%%%%%%%%
\begin{abstract}

The power flow equations are at the core of most of the computations for designing and operating electric power systems. The power flow equations are a system of multivariate nonlinear equations which relate the power injections and voltages in a power system. A plethora of methods have been devised to solve these equations, starting from Newton-based methods to homotopy continuation and other optimization-based methods. While many of these methods often efficiently find a high-voltage, stable solution due to its large basin of attraction, most of the methods struggle to find low-voltage solutions which play significant role in certain stability-related computations. While we do not claim to have exhausted the existing literature on all related methods, this tutorial paper introduces some of the recent advances in methods for solving power flow equations to the wider power systems community as well as bringing attention from the computational mathematics and optimization communities to the power systems problems. After briefly reviewing some of the traditional computational methods used to solve the power flow equations, we focus on three emerging methods: the numerical polynomial homotopy continuation method, Gr\"obner basis techniques, and moment/sum-of-squares relaxations using semidefinite programming. In passing, we also emphasize the importance of an upper bound on the number of solutions of the power flow equations and review the current status of research in this direction.
\end{abstract}

%%%%%%%%%%%%%%%%%%%%%%%%%%%%%%%%%%%%%%%%%%%%%%%%%%%%%%%%%%%%%%%%%%%%%%%%%%%%%%%%

\section{Introduction}
The power flow equations model the steady-state relationship between the complex voltage phasors and the power injections in an electric power system. Power systems typically operate at a ``high-voltage'' solution to the power flow equations that corresponds to a stable equilibrium of the differential-algebraic equations modeling power system dynamics. Section~\ref{l:pf} of this paper gives an overview of the power flow equations.

Calculating a high-voltage power flow solution has been the subject of research efforts for over fifty years. There exist a variety of mature iterative methods, often based on Newton's method and related variants~\cite{tinney1967,stott1974} or the Gauss-Seidel method~\cite{glimn1957}, that are capable of globally solving many practical large-scale power flow problems. Convergence of these iterative methods depends on the selection of an appropriate initialization. Reasonable initializations, such as the solution to a related problem or a ``flat start'' consisting of a $1\angle 0^\circ$ voltage profile, often result in convergence to the high-voltage power flow solution. However, determining an appropriate initialization is challenging when parameters move outside typical operating ranges as may occur with high penetrations of renewable generation and during contingencies. Illustrating this initialization challenge, \cite{thorp1989} and~\cite{thorp1990} demonstrate that the basins of attraction for Newton-based power flow solution methods are fractal in nature.

This challenge has motivated a variety of approaches for improving the robustness of Newton-based methods with respect to the choice of initialization. For instance, \cite{iwamoto1981} calculates an ``optimal multiplier'' at each iteration of the Newton algorithm that prevents divergence. Other approaches consider alternate formulations of the power flow equations, such as~\cite{bromberg2015} which applies a Newton-based iteration to a formulation that includes both voltages and currents.

There has been significant recent interest in alternatives to Newton-based methods. For instance, the Holomorphic Embedding Load Flow method~\cite{trias2012}, which uses analytic continuation theory from complex analysis, is claimed to be capable of reliably finding a stable solution for any feasible set of power flow equations. Approaches based on monotone operator and convex optimization theory~\cite{dj2015a,dj2015b} identify regions with at most one power flow solution. After identifying such a region, any solution contained within can be quickly calculated. Convex relaxation techniques can also calculate power flow solutions~\cite{lavaei_convexpf} and certify infeasibility~\cite{pfcondition,irep2013,hicss2016}. Additionally, progress has been made using ``decoupling'' approximations that facilitate separate analysis of the active and reactive power flow equations~\cite{ilic1992,dorfler2013}.

Many of these approaches focus primarily on calculation of the high-voltage power flow solution. However, other power flow solutions generally exist, often at lower voltages~\cite{tavora1972,klos1975,baillieul1982,guo1994,klos1991}. These solutions are important for many types of power system stability assessments~\cite{venikov1975,tamura1983,ribbens-pavella1985,chiang1987}. There may also exist multiple stable power flow solutions, particularly in the presence of power flow reversal conditions on distribution systems~\cite{turitsyn2014}. Attempts to calculate multiple power flow solutions include the use of a Newton-based algorithm with a range of carefully chosen initializations~\cite{overbye1996,overbye2000}. Using semidefinite relaxations of the power flow equations with objectives that are functions of squared voltage magnitudes, \cite{allerton2011} also identifies multiple power flow solutions. However, these approaches are not guaranteed to find all solutions.

A numerical continuation approach that claims to find all power flow solutions was presented in~\cite{thorp1993}. Since it scales with the number of \emph{actual} rather than \emph{potential} power flow solutions, this approach is computationally tractable for realistic power systems. However, the robustness proof indicating that this approach finds all solutions to all power flow problems is flawed~\cite{chen2011}, and~\cite{counterexample2013} presents a counterexample. This counterexample also invalidates the robustness of the related approach in~\cite{liu2005} for calculating all power flow solutions that have a certain stability property. Recent work~\cite{lesieutre_wu_allerton2015} presents a related method that improves the robustness of~\cite{thorp1993}.

Although not yet scalable to large systems, there are approaches which are guaranteed to find all power flow solutions. As described in Section~\ref{l:nphc}, the most computationally tractable of these methods are based on numerical polynomial homotopy continuation (NPHC). Existing techniques are tractable for systems with up to 14 buses~\cite{mehta2014a,mehta2014b} (and the equivalent of $18$ buses for the related Kuramoto model~\cite{mehta2015algebraic}). The NPHC methods use continuation to trace all complex solutions from a selected ``simple'' polynomial system for which all solutions can be easily calculated to the solutions of the specified target system. In this context, the power flow equations are transformed by splitting real and imaginary parts of the voltage phasors to obtain polynomials in real variables. Thus, only the real solutions to these polynomial equations are physically meaningful. Nevertheless, ensuring recovery of all real solutions requires a number of continuation traces that depends on an upper bound for the number of complex solutions. Apart from being theoretically interesting, tighter upper bounds on the number of power flow solutions would thus improve the computational tractability of NPHC methods. Existing bounds~\cite{baillieul1982,guo1994} are based on calculations of the number of complex solutions for \emph{complete} (i.e., fully connected) networks. Recent work~\cite{acc2016} uses NPHC to find all complex solutions for a variety of small test cases in order to inform conjectured upper bounds that are based on the network topology.

Other methods guaranteed to find all solutions to systems of polynomials (which may therefore be applied to the power flow equations) are the eigenvalue technique in~\cite{dreesen2009}, interval analysis~\cite{mori1999}, Gr\"obner bases (Section~\ref{l:grobner}), and the ``moment/sum-of-squares'' relaxations in~\cite{lasserre_book} and Section~\ref{l:msos}.

Section~\ref{l:grobner} overviews the Gr\"obner basis techniques. In addition to solving the power flow equations, these techniques are presented in the contexts of equivalencing methods, bifurcation analyses, and dynamic power system models.

Section~\ref{l:msos} describes recent progress in moment/sum-of-squares relaxations of the power flow equations. While these relaxations are applicable to many power systems computations, this work is presented in the context of the \emph{optimal power flow} (OPF) problem. Specifically, the OPF problem seeks the voltages which result in power injections that minimize operational cost while satisfying both the power flow equations and engineering limits. The moment/sum-of-squares relaxations lower bound the optimal objective value and, for many OPF problems, give the global solution.

\section{Overview of the Power Flow Equations}
\label{l:pf}

Consider an $n$-bus electric power system where $\mathcal{N} = \left\lbrace 1, \ldots, n \right\rbrace$ is the set of buses and $\mathcal{G}$ is the set of generator buses.\footnote{A ``bus'' in power system terminology represents a node in the graph corresponding to the power system network. A ``line'' corresponds to an edge of this graph.} The network admittance matrix, which contains the electrical parameters of the network as well as the topology information, is denoted $\mathbf{Y} = \mathbf{G} + \mathbf{j} \mathbf{B}$, where $\mathbf{j} = \sqrt{-1}$. (See, e.g.,~\cite{glover_sarma_overbye} for details on the construction of the admittance matrix.)

Each bus has two associated complex values: the voltage phasors and the power injections. We will use both complex voltage phasor representation $V \in \mathbb{C}^n$ and rectangular voltage coordinates $V_{d} + \mathbf{j} V_{q}$, $V_d, V_q \in \mathbb{R}^n$. Each bus~$i \in \mathcal{N}$ has active and reactive power injections $P_{i} + \mathbf{j} Q_{i}$.

In terms of complex voltages, the power flow equations are polynomials in $V$ and $\overline{V}$:
\begin{subequations}
\label{pf_complex}
\begin{align}
\label{pf_P_complex}
& P_i = \Re\left(V_i \sum_{k=1}^n \overline{\mathbf{Y}}_{ik} \overline{V}_k\right) \\
\label{pf_Q_complex}
& Q_i = \Im\left(V_i \sum_{k=1}^n \overline{\mathbf{Y}}_{ik} \overline{V}_k\right)
\end{align}
where $\overline{\left(\,\cdot\,\right)}$ indicates the complex conjugate and $\Re\left(\cdot\right)$ and $\Im\left(\cdot\right)$ return the real and imaginary parts, respectively, of a complex argument. Squared voltage magnitudes are
\begin{align}\label{pf_V_complex}
\left|V_i \right|^2 = V_i \overline{V}_i.
\end{align}
\end{subequations}

Splitting real and imaginary parts of~\eqref{pf_complex} and using rectangular voltage coordinates yields quadratic polynomials in real variables $V_d$ and $V_q$. The active and reactive power injections at bus~$i$ are
\begin{subequations}
\label{pf}
\begin{align}\nonumber
P_i = & V_{di} \sum_{k=1}^n \left( \mathbf{G}_{ik} V_{dk} - \mathbf{B}_{ik} V_{qk} \right) &  &  \\[-3pt] 
\label{pf_P}  & + V_{qi} \sum_{k=1}^n \left( \mathbf{B}_{ik}V_{dk} + \mathbf{G}_{ik}V_{qk} \right) + P_{Di}, \\ \nonumber 
Q_i = & V_{di} \sum_{k=1}^n \left( -\mathbf{B}_{ik}V_{dk} - \mathbf{G}_{ik} V_{qk}\right) \\[-3pt]
\label{pf_Q} & + V_{qi} \sum_{k=1}^n \left( \mathbf{G}_{ik} V_{dk} - \mathbf{B}_{ik} V_{qk}\right) + Q_{Di}.
\end{align}
Squared voltage magnitudes are
\begin{equation} \label{pf_Vsq}
\left|V_i \right|^2 = V_{di}^2 + V_{qi}^2.
\end{equation}
\end{subequations}

To represent typical equipment behavior, each bus is traditionally classified as PQ, PV, or slack. PQ buses, which typically correspond to loads, treat $P_i$ and $Q_i$ as specified quantities and enforce the active and reactive power equations. PV buses, which typically correspond to generators, enforce the active power and squared voltage magnitude equations with specified $P_i$ and $|V_i|^2$. The associated reactive power $Q_i$ may be computed as an ``output quantity'' via the reactive power equation.\footnote{Note that this paper does not consider reactive-power-limited generators.} Finally, a single slack bus is selected with specified $V_i$ (typically chosen such that the reference angle is $0^\circ$; i.e., $\Im\left(V_i\right) = 0$). The active power $P_i$ and reactive power $Q_i$ at the slack bus are determined from the active and reactive power equations, respectively; network-wide conservation of complex power is thereby satisfied. Solving the power flow equations means determining voltage phasors such that the enforced equations are satisfied at each bus.

\section{The Numerical Polynomial Homotopy Continuation Method}
\label{l:nphc}

The Numerical Polynomial Homotopy Continuation (NPHC) method \cite{SVW:96,Li:2003,SW:05} has recently gained attention as it successfully found all the solutions of the power flow equations~\eqref{pf} for the IEEE test systems with up to $14$ buses~\cite{mehta2014a,mehta2014b}. The method has also been applied~\cite{mehta2015algebraic} to find all equilibria of the Kuramoto model~\cite{acebron2005kuramoto}, a prototypical model for the power flow equations \cite{dorfler2011critical}, for up to an $18$-bus case with different network topologies. This section provides an overview of the NPHC method and discusses recent work in determining upper bounds on the number of isolated complex solutions for the power flow equations.

\subsection{Overview of the NPHC Method}
The basic ingredient of the NPHC method, called homotopy continuation, has been used by mathematicians for quite some time~\cite{allgower2003introduction}. To solve a system of nonlinear equations, one first constructs a ``simple'' system for which all solutions can be easily identified. The solutions of the simple system are then continuously deformed to obtain solutions of the original ``target'' system. However, in traditional continuation methods for general nonlinear systems, the solution-paths may cross one another, bifurcate, or turn back towards the simple system. In practice, one or more of these behaviors may be observed for general nonlinear systems. Thus, these methods generally do not guarantee finding all solutions.

For the specific case of systems of \emph{polynomial} equations, however, the situation is different due to the maturity of algebraic geometry theory. For polynomial systems, one can construct a specific type of homotopy method, NPHC, between the simple system and the target system such that all the solution-paths are well-behaved. The NPHC method thus guarantees finding all isolated complex solutions of the target system of polynomial equations. 

We now describe the specific strategy for the NPHC method. Denote the transpose as $\left(\cdot \right)^\intercal$. Let the system of nonlinear polynomial equations $f(x) = 0$, where $f(x)$ is a vector of $n$ equations $f\left(x\right) := \begin{bmatrix}f_1(x) & \ldots & f_{n}(x)\end{bmatrix}^\intercal$ and $x$ is a vector of $n$ variables $x := \begin{bmatrix}x_1 & \ldots & x_n\end{bmatrix}^\intercal$, be the target system to be solved. To utilize theory from algebraic geometry, we consider all the variables $x_i \in \mathbb{C}$. If only real solutions are physically relevant, as is the case for the power flow equations~\eqref{pf}, we disregard all the non-real solutions upon completion of the method. 

The NPHC method starts by calculating an upper bound on the number of isolated complex solutions of \mbox{$f(x) = 0$}. There exist various off-the-shelf upper bounds arising from the computational algebraic geometry literature. The classical B\'ezout bound (CBB) states that the number of complex isolated solutions for a system of polynomial equations is at most the product of degree of each of the polynomials: $\prod_{i=1}^{n} d_{i}$, where $d_i$ is the degree of $f_i(x)$. A discussion on tighter upper bounds for the power flow equations is presented later in this section.

Using an upper bound on the number of complex solutions, one constructs a ``simple'' system in the same variables, \mbox{$g(x) = 0$}, such that 1) the simple system has the same number of isolated complex solutions as the upper bound, and 2) obtaining all the solutions of $g(x) = 0$ is straightforward. Using the CBB, such a system is $g(x) := \begin{bmatrix} a_1 x_1^{d_1} - b_1 & \ldots & a_n x_n^{d_n} - b_n\end{bmatrix}^\intercal$, where $a_i, b_i\neq 0$ are generic complex numbers. The system \mbox{$g\left(x\right) = 0$} has $\prod_{i=1}^{n}d_i$ isolated solutions: $x_i = \sqrt[d_{i}]{\frac{b_{i}}{a_{i}}}$.

The NPHC method then constructs a homotopy between $f(x)$ and $g(x)$:
\begin{equation}
 H(x; t) := (1-t)\, f(x) + \eta\, t\, g(x) = 0,
 \label{eq:PNHC}
\end{equation}
where $t\in [0, 1]$. At $t=1$, we have $g(x) = 0$ and know all the solutions by construction. At $t=0$, we recover the original system $f(x) = 0$. Hence, for each of the solutions of $g(x) = 0$, a path from $t=1$ to $t=0$ is tracked using an efficient path-tracker~\cite{SW:05}, such as a predictor-corrector method. Since the upper bound on the number of solutions of $f(x) = 0$ may not be exact, some of the solution-paths may diverge to infinity along the way. It has been shown that exactly as may paths as the number of isolated complex solutions, counting multiplicity, of $f(x) = 0$ reach $t=0$ so long as $\eta$ in~\eqref{eq:PNHC} is chosen generically from $\mathbb{C}$~\cite{SW:05}. In other words, for $t\in (0,1]$ with a generic complex $\eta$, it is proven that no path will cross other paths, turn back, or bifurcate. (Note that more than one path may sometimes reach the same solution at $t=0$ if the system has multiple solutions.) Hence, in the end, the NPHC method is guaranteed to yield \emph{all} isolated complex solutions of a system of polynomial equations. The method is \emph{embarassingly parallelizeable} since each path can be tracked independently, and is thus suitable for high-performance supercomputing clusters.

\subsection{Upper Bounds on the Number of Power Flow Solutions}
The efficiency of the NPHC strategy strongly depends on the quality of the upper bounds on the isolated complex solutions of the power flow equations, as this is the number of solution-paths one has to track. Another determinant of the efficiency of the method is the computational effort required to compute the upper bound and solve the corresponding start system $g(x) = 0$. The CBB is trivial to compute for any given system, and solving the corresponding start system is also straightforward. However, in practice, the number of complex solutions is usually much lower than the CBB due to the sparsity of the system. Hence, computational effort is wasted in tracking paths that eventually diverge.

To apply NPHC to large systems of power flow equations, we must determine an upper bound on the number of complex solutions that is as tight as possible. One approach for computing high-quality upper bounds is by exploiting the sparsity structure inherent to the power flow equations, which is encoded inside the network topology of the power system.

After fixing the slack bus voltage, there are $2n-1$ power flow equations for an $n$-bus system, each of which is a quadratic polynomial. Thus, the CBB is $2^{2(n-1)}$. A tighter upper bound, $\binom{2n - 2}{n-1}$, on the number of complex solutions was first proposed in~\cite{baillieul1982} for systems without PQ buses and extended to general power systems in~\cite{guo1994}. In~\cite{mehta2014a}, it was pointed out that a generic upper bound, called the BKK bound (named after its inventors Bernstein-Khovanskii-Kushnirenko~\cite{bernshtein1975number,kushnirenko1975newton,khovanskii1978newton}), was significantly smaller than the CBB for the IEEE test case with up to $14$ buses. However, no theoretical justification for this observation is yet known.

These bounds do not consider the sparsity structure of the power flow polynomials. There is only limited work on bounds that are functions of the network topology. The bound in~\cite{guo1990} applies to power systems whose network graphs are composed of subgraphs having exactly one shared bus. The number of isolated complex solutions for such systems is equal to the product of number of solutions for the individual subgraphs. Using numerical experiments conducted with the NPHC method, recent results in~\cite{acc2016} conjecture an upper bound for systems whose network graphs are composed of \emph{maximal cliques} (i.e., maximally sized, completely connected subgraphs) which share exactly two buses along with some other technical conditions. Upper bounds for many other types of network topologies are yet to be discovered. 

\section{The Gr\"obner Basis Techniques}
\label{l:grobner}

As discussed in the previous sections, many important power systems problems are naturally formulated in terms of systems of polynomial equations. The classical and perhaps most important is the system of power flow equations. Other examples include but are not limited to loadability limits, small-signal stability and others. These systems of equations define an algebraic set characterizing the relation between system variables and parameters. Although these relations are usually defined in high-dimensional spaces, the important engineering questions usually involve only a few variables and parameters. For example, classical loadability and path rating analysis addresses the questions of system response to variation of specific load consumption levels. 

Gr\"obner basis techniques provide a formal way to construct all the solutions of polynomial systems equations. Apart from offering an alternative to other techniques described in this tutorial, they also allow elimination of variables from these systems at the expense of raising the total degree of equations. In the following we provide an informal introduction to the algebraic geometry concepts necessary to understand the ideas of Gr\"obner bases. The reader is referred to \cite{cox1992ideals} for an accessible but formal and thorough introduction to the subject.

From algebraic perspective, any polynomial equations in $x_1,\dots, x_n$ variables with real coefficients is an element of the so-called polynomial ring $\mathbb{R}[x_1,\dots, x_n]$ or $\mathbb{R}[\mathbf{x}]$ -- the set of polynomials in $x_1,\dots, x_n$ with real coefficients with the naturally defined addition and multiplication operations. Like elements of other rings, for example integer numbers, the polynomials can be sometimes factorized into products of lower order polynomials, for example $x_1^2-x_2^2 = (x_1 - x_2) (x_1 + x_2)$, which allows the definition of a long division operation on the polynomial ring. Whenever the degree of polynomial $Q \in \mathbb{R}[\mathbf{x}]$ is less than the degree of polynomial $P \in \mathbb{R}[\mathbf{x}]$, one can always represent $P = U Q + R$, where $U$ is the quotient and $R$ the remainder polynomials, each with lower order in comparison to $Q$. As there is no natural ordering on multivariable polynomials, typically a lexicographic ordering is used, where the monomials are ordered lexicographically, so $1 \prec x_1 \prec x_1^2 \prec \dots \prec x_2 \prec  x_2 x_1 \prec x_2 x_1^2 \prec \dots \prec x_2^2$ etc, and polynomials are ordered with respect to the highest monomial order, so, for example $x_1^2 + x_2^2 \succ 2 x_1 x_2$ because $x_2^2 \succ 2 x_1 x_2$. Given a polynomial ordering, the long division is uniquely defined. For example, dividing $P=x_1^2 + x_2^2$ by $Q = x_1 + x_2$ one gets $x_1^2 + x_2^2 = (x_2 - x_1)(x_1 + x_2) + 2 x_1^2$.

An \textit{ideal} $\langle P_1,\dots, P_m\rangle$ generated by polynomials $P_1(\mathbf{x}),\dots, P_m(\mathbf{x}) \in \mathbb{R}[\mathbf{x}]$ is a set of polynomials that can be represented as linear combination of $P_1\dots, P_m$, i.e. $f = g_1 P_1 + \dots + g_m P_m$, where every coefficient $g_k \in \mathbb{R}[\mathbf{x}]$ is itself a polynomial in the same variable set. Any point $\mathbf{y}$ that solves the system of equations $P_1(\mathbf{y}) = 0,\dots, P_m(\mathbf{y}) = 0$ is also a root of every polynomial $P$ in the ideal $\langle P_1,\dots, P_m \rangle$. An ideal defines a set of equations that can be constructed from the original equations $P_k(x)=0$ via algebraic addition and multiplication operations. Any ideal can be generated by various different base polynomials. For a trivial example, consider replacing the first equation $P_1(\mathbf{x})=0$ with $P_1(\mathbf{x}) + P_2(\mathbf{x}) = 0$. Obviously the new set of polynomials will define exactly the same ideal and characterize exactly the same algebraic set of solutions. 

The Gr\"obner basis is a special representation of an ideal that possesses a number of properties making it particularly suitable for solving elimination problems and also for constructing complete sets of solutions. For the lexicographic ordering of monomials defined above, the Gr\"obner basis provides a ``triangular'' representation of the ideal $\langle Q_1, \dots, Q_m\rangle$ with the following properties: $Q_1 = Q_1(x_n)$, $Q_2 = Q_2(x_{n-1}, x_{n})$, and so on with $Q_n = Q_n(x_1\dots, x_n)$ where we assumed for simplicity that $m=n$. Such a representation of the original system of equations provides a straightforward way of finding all the solutions. First, the univariate equation $Q_1(x_n) = 0$ can be solved to find all the values of $x_n$. This can be accomplished, for example, by constructing the companion matrix and finding its eigenvalues, or using the Newton method followed by factorization of the obtained solutions. 

On the second step, the solutions of the first equation can be substituted into the equation $Q_2(x_{n-1}, x_{n}) = 0$ to find all the values of $x_{n-1}$. The obvious advantage of this approach is the transformation of the original multivariate problem to a sequence of univariate polynomial problems for which generally more straightforward methods exist both for identification of the solutions and for providing guarantees that no solution exist in specific regions. Certainly, the degree of the univariate polynomials can become very large for Gr\"obner bases constructed from large systems of equations. For systems of quadratic equations common to power systems, the degrees of the polynomials are bounded by $2^{2^n+1}$ where $n$ is the total number of variables. 

In many practical problems the goal is to characterize the dependence of one or many variables on external parameters. Algebraically, both the unknown variables $x_1, \dots, x_m$ and external parameters $p_1, \dots, p_k$ can be combined in a single variable set $x_1,\dots, x_n$ with $n=m+k$ and $x_{m+i} = p_i$ for $i=1,\dots, k$. In this case the Gr\"obner basis takes the form $Q_1 = Q_1(x_m, x_{m+1},\dots, x_n)$, $Q_2 = Q_2(x_{m-1},\dots, x_n)$, and, finally $Q_m = Q_m(x_1,\dots, x_n)$. The first equations represents the implicit dependence of the unknown variable $x_m$ on external parameters $p_1, \dots, p_k$ represented by $x_{m+1},\dots, x_n$. 

The Gr\"obner basis can be constructed using the classical Buchberger's algorithm which generalizes the classical Euclidean algorithm for finding the greatest common divisor of two integers to polynomial rings. In essence, the Buchberger algorithm implements the nonlinear version of the Gaussian elimination procedure and establishes the analogue of LU decomposition for nonlinear systems of equations. This algorithm is implemented in most of general purpose computer algebra systems like Maple, Mathematica, Sagemath as well as specialized algebraic geometry systems like MAGMA and Singular. The algebraic nature of the algorithm requires exact arithmetic and definition of the problem on the polynomial ring $\mathbb{Q}[\mathbf{x}]$ over field of rational rather than real numbers. This constraint reduces the computational efficiency of the algorithm and limits its applicability to large systems, although at the same time it eliminates the numerical error that can be accumulated in floating point manipulations. Existing implementations of Gr\"obner basis construction algorithms do not exploit any structural properties of the underlying network graphs. However, recent generalization of the classical algorithm suggest that the complexity of the algorithm and degree of the resulting polynomials may be substantially lower than classical upper bounds for graphs with low tree-width typical for power systems \cite{parrilo2014}. So applications of Gr\"obner basis approach to large-scale power systems are not out of question and may become realistic in the future with the maturation of algorithms exploiting the low treewidth of the underlying system graphs.

Below we discuss several potential applications of Gr\"obner basis techniques to power system analysis problems, illustrating the application of Gr\"obner basis techniques on the simple example of a two-bus power system.

\begin{figure}[t]
    \centering{
        \includegraphics[width=\columnwidth]{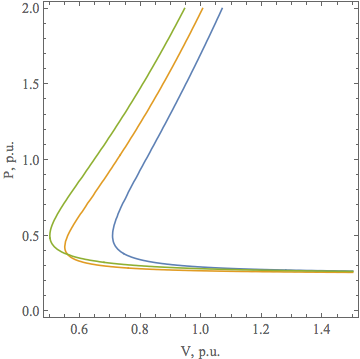}
    }
    \caption{Dependence of the total power consumption $P$ on voltage at bus~$1$ for $P_L = 0.25 \,p.u., Q_L = 0.25 \,p.u.$ (blue),
$P_L = 0.25 \,p.u., Q_L = 0.0 \,p.u.$ (orange), $P_L = 0.25 \,p.u., Q_L = -0.25 \,p.u.$ (green).
    }
    \label{fig:pv}
\end{figure}

The first application we discuss is equivalencing and model reduction. A power grid network has a hierarchical structure with high-voltage power line layers feeding the power to lower-voltage layers. Typically, the low-voltage distribution grid models representing the power consumers are aggregated into simple equivalents typically represented either with constant power $PQ$-buses or slightly more sophisticated ZIP models. The equivalencing of distribution grids happens in an ad-hoc non-systematic manner. Gr\"obner bases provide a tool for constructing the equivalent models in a systematic way. For a trivial example, consider a classical two-bus system where a bus $1$ with voltage  $V$ is connected to a $PQ$ load bus via a line of impedance $r+\mathbf{j} x$. The effective load model as seen by the transmission grid represented by bus $1$ is different from the constant specified $PQ$ load due to losses in the line. The implicit model can be recovered by constructing the Gr\"obner basis for the following system of power flow equations:
\begin{align}\label{eq:twobus-first}
 P & = V i_d \\
 Q & = -V i_q \\
 V_x &= V - r i_d +x i_q \\
 V_y &= -r i_q - x i_d \\
 P_L &= V_d i_d - V_q i_q \\
 Q_L &= -V_d i_q - V_q i_d. \label{eq:twobus-last}
\end{align}
In this system $V_x + \mathbf{j} V_y$ represents voltage phasor at bus~$2$ and $i_x+\mathbf{j}i_y$ represents the current phasor. The power levels $P,Q$ represent the overall power consumption at bus~$1$, while $P_L$, $Q_L$ correspond to consumption at load bus~$2$. The levels of $r$, $x$ are set to $r=x=1/4$~per~unit (p.u.) for the sake of simplicity in this illustration, whereas the values of $V$, $P_L$, and $Q_L$ are assumed to be independent parameters. Using the lexicographic ordering $i_d \prec i_q \prec  V_d \prec V_q \prec Q \prec P \prec V \prec P_L \prec Q_L$, one can construct the Gr\"obner basis for this system and obtain the following first two generating equations:
\begin{align}
& 2 P^2 -2 A P + B = 0  \\
& P - P_L - Q + Q_L = 0,
\end{align}
where $A = P_L - Q_L + 2 V^2$ and $B = P_L^2 - 2 P_L Q_L + Q_L^2 + 4 P_L V^2$. These relations define an implicit model of the aggregate load, representing the dependence of the power consumption levels $P$ and $Q$ on terminal voltage $V$. Note that typically the value of $V$ is affected by the flows in the transmission grid, while external parameters $P_L$ and $Q_L$ depend on the load consumption levels. This model can either be used directly or approximated by more common polynomial models of the loads. The resulting dependence is illustrated in Fig.~\ref{fig:pv}. Note that the solution does not exist for low-voltage levels, which corresponds to phenomena referred to as voltage collapse in power system literature. Also, there are two possible values of power consumption for given voltage levels, corresponding to two branches of power flow equation solutions. In this representation only the lower branch is stable. 

The method can be obviously generalized to more complicated situations with less trivial grid topologies and more sophisticated individual load models. As discussed before, the traditional Gr\"obner basis construction algorithms may result in high-order implicit equations describing the models, so either numerical approximation of intermediate results or utilization of more advanced sparsity-exploiting elimination schemes is necessary to scale the approach to large systems. 

Another important application of the Gr\"obner basis approach is the analysis of bifurcations in power system models. Overloading the power system beyond acceptable limits typically leads to disappearance of equilibrium or loss of stability. Understanding the loadability limits is critical for secure operation of power systems. Mathematically, disappearance of power flow solutions happens through the saddle-node bifurcation which occurs when the power flow equations' Jacobian matrix becomes singular. As the power flow equations have a polynomial form, the Jacobian can be also represented as matrix polynomial $J(x,p)$ depending both on the variables and external parameters. Given this representation, one can introduce the singularity condition through two additional equations $J z = 0$ and $z^\intercal z = 1$ where $z$ is the zero eigenvector of the Jacobian. The combination of power flow equations and these two singularity conditions describes the algebraic saddle-node bifurcation manifold.

\begin{figure}[t]
    \centering{
        \includegraphics[width=\columnwidth]{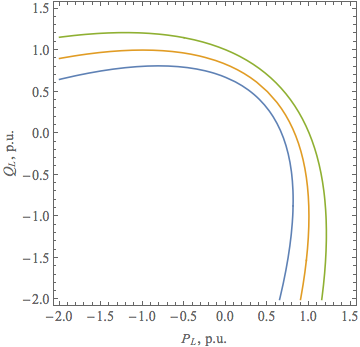}
    }
    \caption{Lodability limit boundaries defined by equation \eqref{eq:loadability} with $V=0.9\, p.u.$ (blue), $V=1.0\, p.u.$ (orange), $V=1.1\, p.u.$ (green).}
    \label{fig:pq}
\end{figure}

For the two-bus example described above, the system of equations \eqref{eq:twobus-first}--\eqref{eq:twobus-last} is complemented with the following representation of the matrix singularity conditions: 
\begin{align}
& z_1 = 0, \\
& z_2 = 0, \\
& z_3 - z_5 i_d+z_6 i_q = 0,\\ 
& z_4 - z_6 i_d-z_5 i_q = 0, \\
& \frac{z_3}{4}+\frac{z_4}{4} - z_5 V_d - z_6 V_q - V z_1= 0, \\
& \frac{z_4}{4}-\frac{z_3}{4} + z_6 V_d - z_5 V_q + V z_2 = 0, \\
& z_1^2+z_2^2+z_3^2+z_4^2+z_5^2+z_6^2 = 1.
\end{align}
The first equation of Gr\"obner basis with the ordering $i_d \prec i_q \prec  V_d \prec V_q \prec Q \prec P \prec z_1 \prec \dots \prec z_6 \prec V \prec P_L \prec Q_L$ takes the form
\begin{equation}\label{eq:loadability}
P_L^2 -2 P_L Q_L + Q_L^2 + 4 V^2 P_L+4 V^2 Q_L - 4 V^4 = 0,
\end{equation}
that describes the ellipsoidal curve corresponding to loadability limits in terms of $P_L$, $Q_L$, and terminal voltage $V$. The resulting boundaries are presented in Fig.~\ref{fig:pq}

Similar techniques can be extended to more complicated problems. For instance, inequality constraints responsible for feasibility can be introduced in the same framework via slack variable approaches. For example, the voltage constraint $V_d^2 + V_q^2 \geqslant \left(V^{min}\right)^2$ can be rewritten in equality form with additional slack variable $s$ introduced as $V_d^2 + V_q^2 = \left(V^{min}\right)^2 + s^2$. This equation has a solution only when the voltage level satisfies the original inequality constraint. Therefore, the bifurcations of this extended system will  define the boundary of the feasibility region where the solution exists and also satisfies the physical limits. 

Gr\"obner basis approaches can be also applied to dynamic models of power systems. For small-signal stability analysis, the linearized equations of motion $\dot x = A(p) x$ depend on the operating conditions described by the parameter set $p$. The small-signal stability condition can be expressed as $A(p) (\mathbf{u} + \mathbf{j} \mathbf{v} ) = (\mu + \mathbf{j} \nu)  (\mathbf{u} + \mathbf{j} \mathbf{v} )$ with extra constraints $u^\intercal u + v^\intercal v = 1$ and $\mu \leqslant 0$. Here, the vector $(u + \mathbf{j} v )$ is the eigenvector of the linearized system matrix, and $\mu,\nu$ are the real and imaginary parts of the corresponding eigenvalue. The condition $\nu \leqslant 0$, which can be also expressed as $\nu + s^2 = 0$ after the slack variable introduction trick, expresses the small-signal stability criterion. The resulting system defines an algebraic manifold corresponding to stable operating conditions. Elimination of the $\mathbf{u},\mathbf{v}$ and $\mu, \nu$ variables defines the manifold in parameter $\mathbf{p}$ space where the system possess a stable equilibrium.
% When expressed in real polynomial equations form this system admits Gr\"obner basis type analysis which allows direct construction of stable regime regions in parameter space. 

Recently, Gr\"obner basis techniques have been applied to power system models in the context of the power flow reversal problem~\cite{turitsyn2014, Nguyen:2015bx}. Introduction of renewable generation on a distribution grid may result in reversal of power flow and export of power to the transmission systems. These new regimes, which are not common to traditional power systems, are characterized by multiple branches of solutions. Although usually not suitable for normal operations, these new equilibrium points may compromise the normal post-fault restoration of the system and result in the system being trapped in an undesirable stable equilibrium. The structure of the new solutions, precise conditions of their appearance, and their impact on power system stability and security is still poorly understood and requires more intensive investigation.

\section{Moment/Sum-of-Squares Relaxations of the Power Flow Equations}
\label{l:msos}

%The power flow equations model the steady-state relationship between the voltages and power injections in a power system and are therefore at the heart of many power system analyses. 

Advances in computational methods related to the power flow equations have the potential to improve solution techniques for many different problems. Motivated by the extensive optimization needs envisioned for future power systems, a variety of convex relaxations and approximations of the power flow equations have recently been developed. These relaxations and approximations can facilitate the optimization of broader classes of systems and operating conditions than traditional methods. Further, emerging computational tools have capabilities that surpass those in traditional solution techniques, such as providing a measure of solution quality, certifying problem infeasibility, and, in many cases, provably obtaining the global optimum. This section reviews recent work in developing a hierarchy of \emph{moment/sum-of-squares} (MSOS) relaxations of the power flow equations.

Most of the work on convex relaxations of the power flow equations has focused on the optimal power flow (OPF) problem (e.g., \cite{lavaei_tps,lesieutre_molzahn_borden_demarco-allerton2011,jabr2011,molzahn_holzer_lesieutre_demarco-large_scale_sdp_opf,low_tutorial,emiliano,andersen2014,coffrin2015,sun2015}). The OPF problem determines a minimum cost operating point for an electric power system subject to both network constraints (i.e., the power flow equations) and engineering limits (e.g., bounds on voltage magnitudes, active and reactive power injections, and line flows). This section presents the MSOS relaxations in terms of the OPF problem; however, note that the MSOS formulations could be applied to a variety of other power system optimization problems (e.g., transmission expansion planning~\cite{jabr2013,taylor2013}, voltage regulation~\cite{zhang2015}, state estimation~\cite{zhu2014}, and calculating voltage stability margins~\cite{pfcondition,irep2013,hicss2016}). Further, when applied to an optimization problem with a feasible space defined by polynomial equality constraints and a constant objective function, MSOS relaxations with sufficiently high relaxation order can find all solutions to the polynomials~\cite{lasserre_book}. This approach could be applied to find all power flow solutions for small systems.

Formulating the power flow equations as a system of polynomials enables the application of the Lasserre hierarchy for polynomial optimization problems~\cite{lasserre_book}. The first-order MSOS relaxation in the Lasserre hierarchy is equivalent to the semidefinite programming (SDP) relaxation in~\cite{lavaei_tps}; higher-order MSOS relaxations in the Lasserre hierarchy take the form of SDPs that generalize the relaxation in~\cite{lavaei_tps}.

Directly applying the Lasserre hierarchy to the OPF problem has been explored in~\cite{pscc2014,cedric,ibm_paper}. Related work includes exploiting sparsity and selectively applying the higher-order constraints to solve large-scale problems~\cite{molzahn_hiskens-sparse_moment_opf,cdc2015}. Other related work develops hierarchies that leverage the \emph{complex structure} of the power flow equations~\cite{complex_hierarchy} and employ a mix of semidefinite and second-order cone programming (SOCP) to enforce the higher-order constraints in the MSOS relaxations~\cite{powertech2015}. This section will review the development of the MSOS relaxations from the Lasserre hierarchy in~\cite{pscc2014,cedric,ibm_paper} and summarize the related work on computational improvements in~\cite{molzahn_hiskens-sparse_moment_opf,cdc2015,complex_hierarchy,powertech2015}.

\subsection{Overview of the Optimal Power Flow Problem}

We first present the non-convex optimal power flow problem, starting with the introduction of notation in addition to that used for the power flow problem in Section~\ref{l:pf}. Superscripts ``max'' and ``min'' denote specified upper and lower limits and buses without generators have maximum and minimum generation set to zero.

In terms of complex voltages, the active and reactive power injections $f_{Pi}\left(V\right)$ and $f_{Qi}\left(V\right)$, respectively, are polynomials in $V$ and $\overline{V}$:
\begin{subequations}
\label{opf_balance_complex}
\begin{align}
& f_{Pi}\left(V\right) := P_{Di} + \Re\left(V_i \sum_{k=1}^n \overline{\mathbf{Y}}_{ik} \overline{V}_k\right) \\
& f_{Qi}\left(V\right) := Q_{Di} + \Im\left(V_i \sum_{k=1}^n \overline{\mathbf{Y}}_{ik} \overline{V}_k\right)
\end{align}
\end{subequations}
where $P_{Di} + \mathbf{j} Q_{Di}$ are specified load demands. Squared voltage magnitudes are
\begin{align}\label{opf_Vsq_complex}
f_{Vi}\left(V\right) := V_i \overline{V}_i.
\end{align}
The OPF problem in terms of complex voltages $V$ is
\begin{subequations}
\label{opf_complex}
\begin{align}
\label{opf_obj_complex} & \min_{V}\quad \sum_{i \in \mathcal{G}} c_i\, f_{Pi}\left(V\right) \qquad \mathrm{subject\; to} \hspace{-160pt} & \\
\label{opf_P_complex} &  \quad P_{Gi}^{\mathrm{min}} \leqslant f_{Pi}\left(V\right) \leqslant P_{Gi}^{\mathrm{max}} & \forall i \in \mathcal{N} \\
\label{opf_Q_complex} &  \quad Q_{Gi}^{\mathrm{min}} \leqslant f_{Qi}\left(V\right) \leqslant Q_{Gi}^{\mathrm{max}} &  \forall i \in \mathcal{N} \\
\label{opf_V_complex} &  \quad (V_{i}^{\mathrm{min}})^2 \leqslant f_{Vi}\left(V\right) \leqslant (V_{i}^{\mathrm{max}})^2 &  \forall i \in \mathcal{N}% \\
%\label{opf_Vref_complex} & \quad \Im\left(V_1\right) = 0
\end{align}
\end{subequations}
where $c \in \mathbb{R}^n$ is a linear cost of active power generation.\footnote{The MSOS relaxations are applicable to more general OPF formulations that include, e.g., line-flow limits, quadratic and piecewise-linear generator costs, and multiple generators per bus~\cite{molzahn_holzer_lesieutre_demarco-large_scale_sdp_opf,andersen2014,molzahn_hiskens-sparse_moment_opf}}

Splitting real and imaginary parts of~\eqref{opf_balance_complex}--\eqref{opf_Vsq_complex} and using rectangular voltage coordinates yields quadratic polynomials in real variables $V_d$ and $V_q$:
\begin{subequations}
\label{opf_balance}
\begin{align}\nonumber
g_{Pi}\left(V_d,V_q\right) & := V_{di} \sum_{k=1}^n \left( \mathbf{G}_{ik} V_{dk} - \mathbf{B}_{ik} V_{qk} \right) &  &  \\[-3pt] 
\label{opf_Pbalance}  & + V_{qi} \sum_{k=1}^n \left( \mathbf{B}_{ik}V_{dk} + \mathbf{G}_{ik}V_{qk} \right) + P_{Di}, \\ \nonumber 
g_{Qi}\left(V_d,V_q\right) & := V_{di} \sum_{k=1}^n \left( -\mathbf{B}_{ik}V_{dk} - \mathbf{G}_{ik} V_{qk}\right) \\[-3pt]
\label{opf_Qbalance} & + V_{qi} \sum_{k=1}^n \left( \mathbf{G}_{ik} V_{dk} - \mathbf{B}_{ik} V_{qk}\right) + Q_{Di}.
\end{align}
\end{subequations}
Squared voltage magnitudes are given by
\begin{equation} \label{opf_Vsq}
g_{Vi}\left(V_d, V_q\right) := V_{di}^2 + V_{qi}^2.
\end{equation}

The OPF problem in terms of real voltage components $V_d$ and $V_q$ is
\begin{subequations}
\label{opf}
\begin{align}
\label{opf_obj} & \min_{V_d,V_q}\quad \sum_{i \in \mathcal{G}} c_i\, g_{Pi}\left(V_d,V_q\right) \qquad \mathrm{subject\; to} \hspace{-160pt} & \\
\label{opf_P} &  \quad P_{Gi}^{\mathrm{min}} \leqslant g_{Pi}\left(V_d,V_q\right) \leqslant P_{Gi}^{\mathrm{max}} & \forall i \in \mathcal{N} \\
\label{opf_Q} &  \quad Q_{Gi}^{\mathrm{min}} \leqslant g_{Qi}\left(V_d,V_q\right) \leqslant Q_{Gi}^{\mathrm{max}} &  \forall i \in \mathcal{N} \\
\label{opf_V} &  \quad (V_{i}^{\mathrm{min}})^2 \leqslant g_{Vi}\left(V_d,V_q\right) \leqslant (V_{i}^{\mathrm{max}})^2 &  \forall i \in \mathcal{N} %\\
%\label{opf_Vref} & \quad V_{q1} = 0.
\end{align}
\end{subequations}

Since the power flow equations are functions of angle \emph{differences}, both formulations of the OPF problem~\eqref{opf_complex} and~\eqref{opf} have a degeneracy in the voltage angle. Similar to the slack bus in the power flow equations, this degeneracy is accounted for by arbitrarily fixing the angle reference at a specified bus. This can either be accomplished by rotating the voltage vectors from the solutions to the MSOS relaxations such that the angle reference is satisfied (as is done here) or by enforcing the additional constraints $\frac{V_1 - \overline{V}_1}{2\mathbf{j}} = 0$ and $V_{q1} = 0$ in~\eqref{opf_complex} and~\eqref{opf}, respectively, where it is assumed that bus~1 has a fixed angle reference of $0^\circ$.

\subsection{Applying the Lasserre Hierarchy to the OPF Problem}
\label{l:msosr}

Lasserre developed a hierarchy of increasingly tighter MSOS relaxations of polynomial optimization problems in real variables~\cite{lasserre_book}.\footnote{The primal form of these relaxations, which is presented here, is derived using truncated \emph{moment} sequences. The dual form is interpreted as a \emph{sum-of-squares} program. Hence, the terminology MSOS relaxations. There is zero duality gap between the primal and dual forms of the MSOS relaxations for OPF problems~\cite{josz-2015}.} MSOS relaxations from the Lasserre hierarchy, which take the form of semidefinite programs, converge to the global optimum of a polynomial optimization problem with increasing relaxation order. We denote the order-$\gamma$ relaxation applied to the OPF problem in real variables $V_d$ and $V_q$ as \msosr{\gamma}. This section reviews~\cite{pscc2014,cedric,ibm_paper}, which develop MSOS relaxations for the OPF problem~\eqref{opf} in real variables.

Development of the relaxations require several definitions. Group the decision variables $V_d$ and $V_q$ into a vector $\hat{x}$:
\begin{equation}
\hat{x} := \begin{bmatrix} V_{d1} & \ldots & V_{dn} & V_{q1} & \ldots & V_{qn}\end{bmatrix}^\intercal.
\end{equation}
Define a vector $x_\gamma$ consisting of all monomials of the voltage components $V_d$ and $V_q$ up to the relaxation order $\gamma$:
\begin{align}
\nonumber
x_\gamma := & \left[ \begin{array}{ccccccc} 1 & V_{d1} & \ldots & V_{qn} & V_{d1}^2 & V_{d1}V_{d2} & \ldots \end{array} \right. \\ \label{x_d}
& \qquad \left.\begin{array}{cccccc} \ldots & V_{qn}^2 & V_{d1}^3 & V_{d1}^2 V_{d2} & \ldots & V_{qn}^\gamma \end{array}\right]^\intercal.
\end{align}

A monomial is defined using a vector $\alpha \in \mathbb{N}^{2n}$ of exponents: $\hat{x}^\alpha := V_{d1}^{\alpha_1} V_{d2}^{\alpha_2}\cdots V_{qn}^{\alpha_{2n}}$. A polynomial is $h\left(\hat{x}\right) := \sum_{\alpha \in \mathbb{N}^{2n}} h_{\alpha} \hat{x}^{\alpha}$, where $h_{\alpha}$ is the real scalar coefficient corresponding to the monomial $\hat{x}^{\alpha}$.

Define a linear functional $L_y\left\lbrace h\right\rbrace$ which replaces the monomials $\hat{x}^{\alpha}$ in a polynomial $h\left(\hat{x}\right)$ with real scalar variables $y$:
\begin{equation}
\label{eq:Lreal}
L_y\left\lbrace h \right\rbrace := \sum_{\alpha \in \mathbb{N}^{2n}} g_{\alpha} y_{\alpha}.
\end{equation}
\noindent For a matrix $h\left(\hat{x}\right)$, $L_y\left\lbrace h\right\rbrace$ is applied componentwise to each element of $h\left(\hat{x}\right)$. 

Consider, e.g., the vector \mbox{$\hat{x} = \begin{bmatrix}V_{d1} & V_{d2} & V_{q1} & V_{q2} \end{bmatrix}^\intercal$} corresponding to the voltage components of a two-bus system. Consider also the polynomial $\left(V_2^{\max}\right)^2 - g_{V2}\left(\hat{x}\right) = \left(V_2^{\max}\right)^2 - V_{d2}^2 - V_{q2}^2$. (The constraint $\left(V_2^{\max}\right)^2 - g_{V2}\left(\hat{x}\right) \geqslant 0$ forces the voltage magnitude at bus~2 to be less than or equal to $V_2^{\max}$.) Then $L_y\left\lbrace \left(V_2^{\max}\right)^2 - g_{V2}\right\rbrace = \left(V_2^{\max}\right)^2y_{0000} - y_{0200} - y_{0002}$. Thus, $L_y\left\lbrace g \right\rbrace$ converts a polynomial $\left(V_2^{\max}\right)^2 - g_{V2}\left(\hat{x}\right)$ to a linear function of $y$.

The MSOS relaxations add additional variables and constraints that are redundant in the original problem but serve to strengthen the relaxation. To simplify the notation, we group all the constraints in~\eqref{opf} into a vector $g\left(\hat{x}\right) \in \mathbb{R}^{6n}$:

\begin{align}\label{ggroup}
& g\left(\hat{x}\right) := \begin{bmatrix} 
P_{G1}^{\max} - g_{P1} \\
\vdots \\
P_{Gn}^{\max} - g_{Pn}\\
g_{P1} - P_{G1}^{\min} \\ 
\vdots \\
g_{Pn} - P_{Gn}^{\min} \\ 
Q_{G1}^{\max} - g_{Q1} \\
\vdots \\
Q_{Gn}^{\max} - g_{Qn}\\
g_{Q1} - Q_{G1}^{\min} \\ 
\vdots \\
g_{Qn} - Q_{Gn}^{\min} \\
\left(V_{1}^{\max}\right)^2 - g_{V1} \\
\vdots \\
\left(V_{n}^{\max}\right)^2 - g_{Vn}\\
g_{V1} - \left(V_{1}^{\min}\right)^2 \\ 
\vdots \\
g_{Vn} - \left(V_{n}^{\min}\right)^2 %\\
%V_{q1} \\
%-V_{q1}
\end{bmatrix}.
\end{align}

The order-$\gamma$ relaxation \msosr{\gamma} can be derived from a rank-constrained optimization problem related to~\eqref{opf}:
\begin{subequations}
\label{opf_aug}
\begin{align}
\label{opf_aug_obj} & \min_{\hat{x}}\quad \sum_{k \in \mathcal{G}} c_k\, g_{Pk}\left(\hat{x}\right) \qquad \mathrm{subject\; to} \hspace{-200pt} & \\
\label{opf_aug_constraint} & \quad g_i\left(\hat{x}\right) x_{\gamma-1}^{\vphantom{\intercal}} x_{\gamma-1}^\intercal \succcurlyeq 0 & i=1,\ldots,6n \\ 
\label{opf_aug_moment}
&\quad x_{\gamma}^{\vphantom{\intercal}} x_{\gamma}^\intercal \succcurlyeq 0
\end{align}
\end{subequations}
where $\succcurlyeq$ indicates positive semidefiniteness of the corresponding matrix and $\gamma \geqslant 1$ is a specified integer which will denote the relaxation order. The rank-one matrices $x_{\gamma-1}^{\vphantom{\intercal}} x_{\gamma-1}^\intercal$ and $x_{\gamma}^{\vphantom{\intercal}} x_{\gamma}^\intercal$ are positive semidefinite by construction. Thus, the constraints $g_i\left(\hat{x}\right) x_{\gamma-1}^{\vphantom{\intercal}} x_{\gamma-1}^\intercal \succcurlyeq 0$ and \mbox{$g_i\left(\hat{x}\right) \geqslant 0$} are redundant and~\eqref{opf_aug_moment} is unnecessary. The matrices $g_i\left(\hat{x}\right) x_{\gamma-1}^{\vphantom{\intercal}} x_{\gamma-1}^\intercal \succcurlyeq 0$ and $x_{\gamma}^{\vphantom{\intercal}} x_{\gamma}^\intercal \succcurlyeq 0$ are composed of polynomials and monomials, respectively, which have degree at most $2\gamma$.

The relaxation \msosr{\gamma} is formed by applying the linear functional $L_y\left\lbrace \cdot \right\rbrace$ to the constraints and objective function in~\eqref{opf_aug} in order to obtain a SDP in terms of the \emph{lifted} variables $y$:
\begin{subequations}
\label{msos}
\begin{align}
\label{msos_obj} & \min_{y}\quad \sum_{k \in \mathcal{G}} c_k\, L_y\left\lbrace g_{Pk}\left(\hat{x}\right) \right\rbrace \qquad \mathrm{subject\; to} \hspace{-200pt} & \\
\label{msos_local} & \quad L_y\left\lbrace g_i\left(\hat{x}\right) x_{\gamma-1}^{\vphantom{\intercal}} x_{\gamma-1}^\intercal \right\rbrace \succcurlyeq 0 & i=1,\ldots,6n \\
\label{msos_moment}
&\quad L_y\left\lbrace x_{\gamma}^{\vphantom{\intercal}} x_{\gamma}^\intercal \right\rbrace \succcurlyeq 0 \\
\label{y0}
& \quad y_{0\ldots 0} = 1.
\end{align}
\end{subequations}
In Lasserre's terminology~\cite{lasserre_book},~\eqref{msos_local} constrains \emph{localizing matrices} and~\eqref{msos_moment} constrains the \emph{moment matrix} for \msosr{\gamma}. See~\cite{pscc2014} for the localizing and moment matrices for small example OPF problems. Note that~\eqref{y0} results from the fact that $
\hat{x}^0 = 1$.

Since $x_0$ is the scalar 1, the localizing constraints~\eqref{msos_local} for \msosr{\gamma} are in fact non-negativity constraints $L_y\left\lbrace g_i\left(\hat{x}\right) \right\rbrace \geqslant 0$. Thus, \msosr{1} is equivalent to the SDP relaxation in~\cite{lavaei_tps}, and \msosr{\gamma} generalizes the relaxation of~\cite{lavaei_tps} for $\gamma > 1$. 

As a relaxation, the optimal objective value for \msosr{\gamma} lower bounds the globally optimal objective value for the OPF problem. A solution to \msosr{\gamma} which satisfies 
\begin{equation} \label{rankcondition}
\mathrm{rank}\left(L_y\left\lbrace \hat{x}\hat{x}^\intercal \right\rbrace\right) = 1
\end{equation}
indicates that the relaxation is \emph{exact}. The globally optimal voltage phasor solution $V^\ast$ to the OPF problem can be extracted using a spectral decomposition. Let $\lambda$ and $\eta$ denote the non-zero eigenvalue and corresponding unit-length eigenvector, respectively, of the matrix $L_y\left\lbrace \hat{x}\hat{x}^\intercal \right\rbrace$. Then $V^\ast = \sqrt{\lambda} \left(\eta_{1:n} + \mathbf{j} \eta_{n+1:2n} \right)$, rotated such that $\Im\left(V_1^\ast\right) = 0$, where subscripts indicate vector entries in MATLAB notation. If the rank condition~\eqref{rankcondition} is not satisfied, increasing the relaxation order to $\gamma+1$ will tighten the relaxation and may yield the globally optimal solution. The results in~\cite{pscc2014,cedric,ibm_paper} demonstrate that \msosr{2} is capable of globally solving many small problems for which \msosr{1} (or, equivalently, the SDP relaxation in~\cite{lavaei_tps}) fails to be exact.

\subsection{Computational Improvements}

Application to large OPF problems requires addressing the computational scaling challenges inherent to the MSOS relaxations. For an $n$-bus system, the size of the moment matrix~\eqref{msos_moment} for \msosr{\gamma} is $\left(2n+\gamma\right)! / \left( \left(2n\right)! \gamma!\right)$. For example, $n=10$ and $\gamma=3$ correspond to a moment matrix of size $1,\!771\times 1,\!771$. This limits application of the ``dense'' formulation described in Section~\ref{l:msosr} to systems with at most approximately ten buses. We next summarize three approaches for improving the computational tractability of the MSOS relaxations: 1) exploiting sparsity and selectively applying the higher-order relaxation constraints to specific ``problematic'' buses, 2) using a hierarchy that exploits the complex structure of the OPF problem, and 3) enforcing the first-order constraints with the SDP formulation while using a SOCP relaxation of the higher-order constraints.

\subsubsection{Exploiting Sparsity and Selectively Applying the Higher-Order Constraints}

As is common in power system optimization, exploiting sparsity has the potential to improve the computational tractability of the MSOS relaxations. A method for exploiting sparsity in relaxations of general polynomial optimization problems~\cite{waki2006} can be directly applied to the MSOS relaxations of the OPF problem. This approach uses a matrix completion theorem~\cite{gron1984} to decompose the single large positive semidefinite constraints in~\eqref{msos_local} and~\eqref{msos_moment} into constraints on many smaller matrices in a manner that depends on the so-called \emph{chordal sparsity} of the power system network. Related approaches applied to the SDP relaxation of~\cite{lavaei_tps} (and therefore also \msosr{1}) enable solution of problems with thousands of buses~\cite{jabr2011,molzahn_holzer_lesieutre_demarco-large_scale_sdp_opf}.

Na\"ively applying the chordal-sparsity approach to the higher-order relaxations is less successful; only systems with at most approximately 40 buses are computationally tractable for \msosr{2}. The key insight needed to scale the higher-order MSOS relaxations to larger problem is that the computationally intensive higher-order constraints are only needed for certain buses. By selectively applying the higher-order constraints at ``problematic'' buses, \cite{molzahn_hiskens-sparse_moment_opf} and~\cite{cdc2015} are able to globally solve OPF problems with thousands of buses.\footnote{More so than computational speed, numerical convergence limits the performance of the higher-order relaxations. Reference \cite{cdc2015} presents a method for preprocessing the OPF problem data to eliminate ``low-impedance'' lines in order to improve convergence characteristics of the SDP solver.} Reference~\cite{molzahn_hiskens-sparse_moment_opf} describes a iterative method for identifying the problematic buses using a ``power injection mismatch'' heuristic. See~\cite{pscc2016} for a computational analysis of the key parameter in this sparsity-exploiting approach (i.e., the number of additional buses with higher-order constraints applied at each iteration).

\subsubsection{Complex Hierarchy}
\label{l:msosc}

The approach in Section~\ref{l:msosr} converts the OPF problem from a formulation with complex variables~\eqref{opf_complex} to a formulation with real variables~\eqref{opf} and then applies the Lasserre hierarchy. An alternative approach recognizes that, in general, the operations of converting to the real representation and applying the Lasserre hierarchy are not commutative. For the OPF problem, it is often computationally advantageous to directly build a \emph{complex} moment/sum-of-squares hierarchy, where we denote the order-$\gamma$ relaxation as \msosc{\gamma}, and then convert to real variables before passing the formulation to the SDP solver. This section briefly describes the complex moment/sum-of-squares hierarchy which was first proposed in~\cite{complex_hierarchy}.

This section adopts notation and follows the development of the real moment/sum-of-squares hierarchy in Section~\mbox{\ref{l:msosr}}. A complex monomial is defined using two vectors of exponents $\alpha,\, \beta \in \mathbb{N}^n$: $V^{\alpha} \overline{V}^{\beta} := V_1^{\alpha_1}\cdots V_n^{\alpha_n} \overline{V_1}^{\beta_1} \cdots \overline{V_n}^{\beta_n}$. A polynomial $h\left(V\right) := \sum_{\alpha,\beta \in \mathbb{N}^n} h_{\alpha,\beta} V^{\alpha} \overline{V}^{\beta}$, where $h_{\alpha,\beta}$ is the complex scalar coefficient corresponding to the monomial $V^\alpha \overline{V}^\beta$. Since $h\left(V\right)$ is a real quantity, $\overline{h_{\alpha,\beta}} = h_{\beta,\alpha}$.

Define a linear functional $\hat{L}_{\hat{y}}\left\lbrace h\right\rbrace$ which replaces the monomials in a polynomial $h\left(V\right)$ with complex scalar variables $\hat{y}$:
\begin{equation}
\label{eq:Lcomp}
\hat{L}_{\hat{y}}\left\lbrace h \right\rbrace := \sum_{\alpha,\beta \in \mathbb{N}^{n}} h_{\alpha,\beta}\, \hat{y}_{\alpha,\beta}.
\end{equation}
For a matrix $h\left(V\right)$, $\hat{L}_{\hat{y}}\left\lbrace h\right\rbrace$ is applied componentwise. 

Consider, e.g., the vector $V = \begin{bmatrix}V_{1} & V_{2} \end{bmatrix}^\intercal$ of complex voltages for a two-bus system and the polynomial $\left(V_2^{\max}\right)^2 - f_{V2}\left(V\right) = \left(V_2^{\max}\right)^2 - V_2 \overline{V_{2}}$. (The constraint $\left(V_2^{\max}\right)^2 - f_{V2}\left(V\right) \geqslant 0$ forces the voltage magnitude at bus~2 to be less than or equal to $V_2^{\max}$.) Then $\hat{L}_{\hat{y}}\left\lbrace \left(V_2^{\max}\right)^2 - f_{V2}\left(V\right)\right\rbrace = \left(V_2^{\max}\right)^2\hat{y}_{00,00} - \hat{y}_{01,01}$. Thus, $\hat{L}_{\hat{y}}\left\lbrace h \right\rbrace$ converts a polynomial $h\left(V\right)$ to a linear function of $\hat{y}$.

For the order-$\gamma$ relaxation, define a vector $z_\gamma$ consisting of all monomials of the voltages up to order $\gamma$ without complex conjugate terms (i.e., $\beta = 00\cdots 0$):
\begin{align}
\nonumber
z_\gamma := & \left[ \begin{array}{ccccccc} 1 & V_{1} & \ldots & V_{n} & V_{1}^2 & V_{1}V_{2} & \ldots \end{array} \right. \\ \label{eq:z_d}
& \qquad \left.\begin{array}{cccccc} \ldots & V_{n}^2 & V_{1}^3 & V_{1}^2 V_{2} & \ldots & V_{n}^\gamma \end{array}\right]^\intercal.
\end{align}

We again simplify the notation by grouping all the constraints in~\eqref{opf_complex} into a vector $f\left(V\right) \in \mathbb{R}^{6n}$ in the same manner as in~\eqref{ggroup}.

The order-$\gamma$ relaxation in the complex hierarchy, \msosc{\gamma}, can be derived from a rank-constrained optimization problem related to~\eqref{opf_complex}:
\begin{subequations}
\label{opf_aug_complex}
\begin{align}
\label{opf_aug_obj_complex} & \min_{V}\quad \sum_{k \in \mathcal{G}} c_k\, f_{Pk}\left(V\right) \qquad \mathrm{subject\; to} \hspace{-200pt} & \\
\label{opf_aug_constraint_complex} & \quad f_i\left(V\right) z_{\gamma-1}^{\vphantom{H}} z_{\gamma-1}^H \succcurlyeq 0 & i=1,\ldots,3n \\
\label{opf_aug_moment_complex}
&\quad z_{\gamma}^{\vphantom{\intercal}} z_{\gamma}^H \succcurlyeq 0
\end{align}
\end{subequations}
where $\left(\cdot \right)^H$ denotes the complex conjugate transpose. The rank-one Hermitian matrices $z_{\gamma-1}^{\vphantom{\intercal}} z_{\gamma-1}^H$ and $z_{\gamma}^{\vphantom{H}} z_{\gamma}^H$ are positive semidefinite by construction. Thus, the constraints $f_i\left(V\right) z_{\gamma-1}^{\vphantom{H}} z_{\gamma-1}^H \succcurlyeq 0$ and $f_i\left(V\right) \geqslant 0$ are redundant and~\eqref{opf_aug_moment_complex} is unnecessary. The matrices $f_i\left(V\right) z_{\gamma-1}^{\vphantom{H}} z_{\gamma-1}^H \succcurlyeq 0$ and $x_{\gamma}^{\vphantom{H}} x_{\gamma}^H \succcurlyeq 0$ are composed of complex polynomials and monomials, respectively, which have degree at most $2\gamma$ (i.e., all monomials $V^{\alpha}\overline{V}^{\beta}$ satisfy $\left|\alpha\right| + \left|\beta\right| \leqslant 2\gamma$, where $\left|\,\cdot\, \right|$ indicates the one-norm).

The relaxation \msosc{\gamma} is formed by applying the linear functional $\hat{L}_{\hat{y}}\left\lbrace \cdot \right\rbrace$ to the constraints and objective function in~\eqref{opf_aug_complex} in order to obtain a SDP in terms of the lifted \emph{complex} variables $\hat{y}$:
\begin{subequations}
\label{msos_complex}
\begin{align}
\label{msos_obj_complex} & \min_{\hat{y}}\quad \sum_{k \in \mathcal{G}} c_k\, \hat{L}_{\hat{y}}\left\lbrace f_{Pk}\left(V\right) \right\rbrace \qquad \mathrm{subject\; to} \hspace{-200pt} & \\
\label{msos_local_complex} & \quad \hat{L}_{\hat{y}}\left\lbrace f_i\left(V\right) z_{\gamma-1}^{\vphantom{H}} z_{\gamma-1}^H \right\rbrace \succcurlyeq 0 & i=1,\ldots,3n+2 \\
\label{msos_moment_complex}
&\quad \hat{L}_{\hat{y}}\left\lbrace z_{\gamma}^{\vphantom{H}} z_{\gamma}^H \right\rbrace \succcurlyeq 0 \\
\label{y0_complex}
& \quad \hat{y}_{0\ldots 0,0\ldots 0} = 1.
\end{align}
\end{subequations}
Mirroring Lasserre's terminology~\cite{lasserre_book},~\eqref{msos_local} constrains complex localizing matrices and~\eqref{msos_moment} constrains the complex moment matrix for \msosc{\gamma}. See~\cite{pscc2016} for the complex localizing and moment matrices for small example problems.

Similar to \msosr{\gamma}, a solution to \msosr{\gamma} that satisfies the rank condition
\begin{equation} \label{rankcondition_complex}
\mathrm{rank}\left(\hat{L}_{\hat{y}}\left\lbrace V\, V^H \right\rbrace\right) = 1
\end{equation}
indicates that the relaxation is \emph{exact} and the globally optimal voltage phasor solution $V^\ast$ to the OPF problem can be extracted using a spectral decomposition. Let $\hat{\lambda}$ and $\hat{\eta}$ denote the non-zero eigenvalue and corresponding unit-length eigenvector, respectively, of the matrix $\hat{L}_{\hat{y}}\left\lbrace V\, V^H \right\rbrace$. Then $V^\ast = \sqrt{\hat{\lambda}} \hat{\eta}$, rotated so that $\Im\left(V^\ast_1\right) = 0$. If the rank condition~\eqref{rankcondition_complex} is not satisfied, increasing the relaxation order to $\gamma+1$ will tighten the relaxation and may yield the globally optimal solution. The second-order complex hierarchy globally solves many problems for which the first-order relaxation is not exact.

Since $z_0$ is the scalar 1, the localizing constraints~\eqref{msos_local} for \msosc{\gamma} are in fact non-negativity constraints $\hat{L}_{\hat{y}}\left\lbrace f_i\left(V\right) \right\rbrace \geqslant 0$. It is proven in~\cite{complex_hierarchy} that \msosc{1} and \msosr{1} (and therefore also the SDP relaxation of~\cite{lavaei_tps}) have the same optimal objective values. However, \msosc{1} typically has computational advantages over \msosr{1}~\cite{complex_hierarchy}.

Another characteristic of the complex hierarchy, as proven in~\cite{complex_hierarchy}, is a sum-of-squares interpretation of the dual problem, with zero duality gap between the primal and dual relaxations for the OPF problem.

As proven in~\cite{complex_hierarchy}, the optimal objective value from the real hierarchy is at least as tight as that from the complex hierarchy, and there exist general complex polynomial optimization problems for which the real hierarchy is tighter. However, \cite{complex_hierarchy} conjectures that the real and complex hierarchies are equally tight for a class of problems modeling oscillatory phenomena (which includes the OPF problem) augmented with a redundant ``sphere constraint''.

The computational advantage of the complex hierarchy is most evident when comparing the moment matrix sizes between \msosr{\gamma} and \msosc{\gamma}. For an $n$-bus system, the size of the moment matrix~\eqref{msos_moment_complex} for \msosc{\gamma} (converted to real representation for input to the solver~\cite[Ex. 4.42]{boyd2009}) is $2\left(\left(n+\gamma\right)!\right)/\left(n!\gamma!\right)$, which is significantly smaller than $\left(2n+\gamma\right)! / \left( \left(2n\right)! \gamma!\right)$ for \msosr{\gamma}. For example, $n=10$ and $\gamma=3$ correspond to matrices of size $1,\!771\times 1,\!771$ and $572\times 572$ for the real and complex hierarchies, respectively. Note that the sparsity-exploiting approach applied to \msosr{\gamma}~\cite{molzahn_hiskens-sparse_moment_opf} can be adopted to \msosc{\gamma}~\cite{complex_hierarchy}. The results in~\cite{complex_hierarchy} demonstrate that computational speed improvements for the complex hierarchy vs. the real hierarchy can exceed an order of magnitude for some large test cases.

\subsubsection{Mixed SDP/SOCP Hierarchy}
\label{l:mixed_sdpsocp}

First proposed in~\cite{powertech2015}, a relaxation of the Lasserre hierarchy that mixes SDP and SOCP constraints facilities the development of a related hierarchy with tightness and computational burden between \msosr{1} and \msosr{2}.\footnote{Future work includes extending the mixed SDP/SOCP hierarchy reported in this section to the complex hierarchy \msosc{\gamma}.} This enables exploitation of the intuition that the second-order relaxation is often more than necessary for global solution of many OPF problems.

A necessary (but not sufficient) condition for a matrix to be positive semidefinite takes the form of a SOCP constraint. Specifically, a symmetric, positive semidefinite matrix $\mathbf{W}$ satisfies
\begin{align}\label{socpnecessary}
& \mathbf{W}_{ii} \mathbf{W}_{kk} \geqslant \left|\mathbf{W}_{ik}\right|^2 & \forall \left\lbrace\left(i,k \right) \;|\; k > i \right\rbrace
\end{align}

In~\cite{low_tutorial}, the SOCP constraint~\eqref{socpnecessary} is applied to the first-order relaxation \msosc{1}. While this significantly reduces the computational burden compared to using SDP constraints, the SOCP relaxation in~\cite{low_tutorial} typically only yields the global solution to a limited set of OPF problems.\footnote{The SOCP relaxation in~\cite{low_tutorial} is guaranteed to globally solve OPF problems with radial networks that satisfy certain non-trivial technical conditions, but generally fails for mesh network topologies.} 

Conversely, the mixed SDP/SOCP relaxation proposed in~\cite{powertech2015} formulates the first-order relaxation \msosr{1} with SDP constraints. This alone is sufficient to globally solve many OPF problems~\cite{lavaei_tps,molzahn_holzer_lesieutre_demarco-large_scale_sdp_opf}. When the solution to \msosr{1} does not satisfy the rank condition~\eqref{rankcondition}, the mixed SDP/SOCP relaxation applies the SOCP formulation~\eqref{socpnecessary} to the entries of the localizing and moment matrices which correspond to the higher-order constraints rather than use the computationally intensive SDP formulation. Thus, the proposed mixed SDP/SOCP relaxation forms a ``middle ground'' between the first- and higher-order moment relaxations.

The mixed SDP/SOCP relaxations are generally not as tight as those which use only SDP constraints; there exist OPF problems for which \msosr{2} is exact, but low-order mixed SDP/SOCP relaxations only provide a strict lower bound on the optimal objective value~\cite{illustrative_example}. However, results for a variety of systems with up to several hundred buses demonstrate a significant computational speed improvement for some problems (up to a factor of 18.7 for one test case).

\section{Conclusion}
After briefly reviewing both canonical and emerging computational methods to solve power flow equations, this tutorial paper has presented overviews of three methods based on algebraic geometry: Numerical Polynomial Homotopy Continuation, Gr\"obner Basis, and Moment/Sum-of-Squares relaxations. We anticipate that this tutorial paper will motivate both power systems researchers and computational mathematicians to further explore the capabilities of these methods with respect to power systems problems as well as to enhance the computational capabilities of these methods.

%While discussing the methods, we have also remained explicit in problem formulations so that the reader from other subareas than the power flow equations as well as computational mathematicians can understand and follow the details and intricacies of the problem. We anticipate that this tutorial will serve as a junction between power systems and computational mathematics communities.

%%%%%%%%%%%%%%%%%%%%%%%%%%%%%%%%%%%%%%%%%%%%%%%%%%%%%%%%%%%%%%%%%%%%%%%%%%%%%%%%

\bibliographystyle{IEEEtran}
\bibliography{acc2016_tutorial,bibliography_NPHC_NAG}{}

\end{document}